

\baselineskip=14pt
\parskip=10pt

\magnification=\magstephalf

\def\1{{\overline{1}}}
\def\2{{\overline{2}}}
\parindent=0pt
\overfullrule=0in

\def\frac#1#2{{#1 \over #2}}
\centerline
{\bf 
Explicit Expressions for Moments of the Duration of a 3-Player Gambler's Ruin
}

\bigskip
\centerline
{\it Shalosh B. EKHAD and Doron ZEILBERGER}
\bigskip
\qquad \qquad \qquad {\it Dedicated to Jean-Paul Delahaye, the French analog of Martin Gardener}

{\bf Abstract:} Using {\it experimental mathematics} and {\it symbolic computation}, we derive many moments for the duration of
a $3$-player (fair) gambler's ruin.

{\bf How It All Started: Jean-Paul Delahaye's Wonderful Monthly Column}

We are fortunate to be on the emailing list of Jean-Paul Delahaye, who sends us each month his fascinating {\it Logique et Calcul} column for
the French analog of {\it Scientific American}, called {\it Pour la Science}. The August 2023 piece [De] was dedicated to
the Gambler's Ruin problem, that is one of our favorite topics.
The {\it classical} {\bf two-player} case goes back  to 1657, when Christiaan Huygens  solved it.
It is covered extensively in Feller's {\it probability bible} [F].

Thanks to Delahaye's column, we learned about exciting new research on the {\it three-player case}, by Persi Diaconis and Stewart Ethier [DiE], and, via the latter paper,  about earlier work
[E] [BLT] [S] [R].

{\bf Warm-Up: On the Duration of The 2-player Gambler's Ruin}

Alice and Bob start out with $A$ dollars and $B$ dollars respectively. A fair coin is tossed, and if it lands on Heads, Alice gives Bob a dollar, so now
their capitals are $A-1$ and $B+1$ respectively. On the other hand, if it lands on Tails, Bob gives Alice a dollar, and now
their respective capitals are $A+1$ and $B-1$. Sooner or later, with {\bf probability one}, one of them goes broke, and the
other one goes home with $A+B$ dollars. The duration of the game is a certain {\bf random variable}, let's call is $D$. Huygnes  proved (see [F]) the
famous result
$$
E[D] \,= \, AB \quad.
$$
Less well known is an explicit expression for the {\it variance} (see [F])
$$
Var(D)= \frac{A B \left(A^{2}+B^{2}-2\right)}{3} \quad .
$$

In [Z] we computed the first ten moments (about the mean), but now we have forty(!) of them, and can easily compute more. Let's state a few of them here.
Let $m_r$ be the $r$-th moment about the mean. 
We have
\vfill\eject
$$
m_3(D)=\frac{A B \left(3 A^{4}+10 A^{2} B^{2}+3 B^{4}-10 A^{2}-10 B^{2}+4\right)}{15} \quad ,
$$
\bigskip
$$
m_4(D)= \frac{AB}{105} \cdot
$$
$$
(17 A^{6}+35 A^{5} B +119 A^{4} B^{2}+70 A^{3} B^{3}+119 A^{2} B^{4}+35 A \,B^{5}+17 B^{6}-84 A^{4}-140 A^{3} B -280 A^{2} B^{2}
$$
$$
-140 A \,B^{3}-84 B^{4}+84 A^{2}+140 A B +84 B^{2}+8) \quad,
$$
\bigskip
$$
m_5(D)=\frac{AB}{945} \cdot
$$
$$
(155 A^{8}+630 A^{7} B +1860 A^{6} B^{2}+2730 A^{5} B^{3}+3906 A^{4} B^{4}+2730 A^{3} B^{5}+1860 A^{2} B^{6}+630 A \,B^{7}+155 B^{8}
$$
$$
-1020 A^{6}-3360 A^{5} B -7140 A^{4} B^{2}-8400 A^{3} B^{3}-7140 A^{2} B^{4}-3360 A \,B^{5}-1020 B^{6}+1764 A^{4}
$$
$$
+5040 A^{3} B +5880 A^{2} B^{2}+5040 A \,B^{3}+1764 B^{4}-440 A^{2}-1680 A B -440 B^{2}-144) \quad ,
$$
\bigskip
$$
m_6(D)=\frac{AB}{10395} \cdot
$$
$$
(2073 A^{10}+12573 A^{9} B +43780 A^{8} B^{2}+95040 A^{7} B^{3}+154143 A^{6} B^{4}+172326 A^{5} B^{5}+154143 A^{4} B^{6}+95040 A^{3} B^{7}
$$
$$
+43780 A^{2} B^{8}+12573 A \,B^{9}+2073 B^{10}-17050 A^{8}-86130 A^{7} B -239250 A^{6} B^{2}-418110 A^{5} B^{3}-498960 A^{4} B^{4}-
$$
$$
418110 A^{3} B^{5}-239250 A^{2} B^{6}-86130 A \,B^{7}-17050 B^{8}+42636 A^{6}+182028 A^{5} B 
$$
$$
+367752 A^{4} B^{2}+489720 A^{3} B^{3}+367752 A^{2} B^{4}+182028 A \,B^{5}+42636 B^{6}-29040 A^{4}-116160 A^{3} B 
$$
$$
-143000 A^{2} B^{2}-116160 A \,B^{3}-29040 B^{4}-3344 A^{2}-528 A B -3344 B^{2}-1440) \quad .
$$

To see all the moments up to the $40$-th, look at the output file

{\tt https://sites.math.rutgers.edu/\~{}zeilberg/tokhniot/oGR3moms3a.txt} \quad .

{\bf Arthur Engel's Beautiful Expression for the Expected Duration of a (fair) $3$-player Gambler's Ruin}

In a gambler's ruin game with {\bf three} players, the players have {\it initial capitals} $A,B,C$.
At each round, one of the players is chosen {\it uniformly at random} (i.e. with probability $\frac{1}{3}$) to
give one dollar to one of the other two players (each of whom is equally likely to be the recipient).
Sooner or later, one of the players goes broke. How soon?

Equivalently (and that is how Engel thought about it) you have three towers with $A$, $B$, and $C$ pegs, and at each round, uniformly at random,
you remove a peg from one of the towers and put it (each with probability $\frac{1}{2}$), on one of the other towers
(It is like the Tower of Hanoi, but all pegs have the same diameter). Sooner or later, one of the towers is empty.
How soon?

Using {\it simulation}, and {\it guessing}, Arthur Engel [E] conjectured that the {\it expected duration}, $E[D]$, is given by the
lovely, explicit formula
$$
E[D] \, = \, \frac{3ABC}{A+B+C} \quad .
$$

Once discovered, he had no trouble {\bf proving} it, since calling this quantity $E(A,B,C)$, it obviously satisfies the equation

$$
E(A,B,C) \, = \, 1 \, + \, \frac{1}{6} \, ( E(A+1,B-1,C) \,+ \, E(A-1,B+1,C)
$$
$$
+\,E(A+1,B,C-1) \,+ \,E(A-1,B,C+1)
$$
$$
+\,E(A,B+1,C-1)\, + \,E(A,B-1,C+1) ) \quad ,
$$
subject to the {\bf boundary conditions}
$$
E(0,B,C)=0 \quad, \quad E(A,0,C)=0 \quad, \quad E(A,B,0)=0  \quad,
$$
and this is a {\bf routine verification}. Proving {\bf uniqueness} is fairly easy, using a {\it maximum principle} (see [A], p. 40).

But, in {\bf hindsight}, there was no need for simulation. Because of the boundary conditions, it is natural to {\it postulate} the {\it ansatz}
$$
E(A,B,C) \, = \,\frac{ k ABC}{A+B+C} \quad ,
$$
(note that the denominator $A+B+C$ is the total number of pegs (dollars) and is an invariant). 
Here $k$ is an {\it undetermined} coefficient, that is {\bf to be determined}. Plugging it in the equation, we get that indeed $k=3$.

Bruss et. al. [BLT] worked pretty hard to find the {\it variance}. 
$$
Var(D)= \frac{3 A B C \left(A^{2} B +A^{2} C +A \,B^{2}-3 A B C +A \,C^{2}+B^{2} C +B \,C^{2}-A -B -C \right)}{2 \left(A +B +C \right)^{2}} \quad .
$$

David Stirzaker [S] later gave a simpler proof, using martingales. We will soon describe an even simpler (and much more elementary!) proof, using computer algebra. Using our
method, we will also discover, and prove, explicit expressions for many higher moments (and one can keep going indefinitely).

But let's first state some of these beautiful formulas.

$$
m_3(D)= \frac{3ABC}{40 (A+B+C)^3} \cdot
$$
$$
(9 A^{5} B +9 A^{5} C +54 A^{4} B^{2}+108 A^{4} B C +54 A^{4} C^{2}+90 A^{3} B^{3}-225 A^{3} B^{2} C -225 A^{3} B \,C^{2}+90 A^{3} C^{3}+54 A^{2} B^{4}
$$
$$
-225 A^{2} B^{3} C +1080 A^{2} B^{2} C^{2}-225 A^{2} B \,C^{3}+54 A^{2} C^{4}+9 A \,B^{5}+108 A \,B^{4} C -225 A \,B^{3} C^{2}-225 A \,B^{2} C^{3}+108 A B \,C^{4}
$$
$$
+9 A \,C^{5}+9 B^{5} C +54 B^{4} C^{2} +90 B^{3} C^{3}+54 B^{2} C^{4}+9 B \,C^{5}-90 A^{3} B -90 A^{3} C -180 A^{2} B^{2}+90 A^{2} B C -180 A^{2} C^{2}
$$
$$
-90 A \,B^{3}+90 A \,B^{2} C +90 A B \,C^{2}-90 A \,C^{3}-90 B^{3} C -180 B^{2} C^{2}-90 B \,C^{3}+12 A^{2}
$$
$$
+24 A B +24 A C +12 B^{2}+24 B C +12 C^{2}) \quad ,
$$
\bigskip
$$
m_4(D)= \frac{3ABC}{280 (A+B+C)^4} \cdot
$$
$$
(6 A^{8} B +6 A^{8} C +54 A^{7} B^{2}+108 A^{7} B C +54 A^{7} C^{2}+195 A^{6} B^{3}+354 A^{6} B^{2} C +354 A^{6} B \,C^{2}+195 A^{6} C^{3}
$$
$$
+357 A^{5} B^{4}+63 A^{5} B^{3} C -588 A^{5} B^{2} C^{2}+63 A^{5} B \,C^{3}+357 A^{5} C^{4}+357 A^{4} B^{5}-378 A^{4} B^{4} C 
$$
$$
+3423 A^{4} B^{3} C^{2}+3423 A^{4} B^{2} C^{3}-378 A^{4} B \,C^{4}+357 A^{4} C^{5}+195 A^{3} B^{6}+63 A^{3} B^{5} C +3423 A^{3} B^{4} C^{2}
$$
$$
-8505 A^{3} B^{3} C^{3}+3423 A^{3} B^{2} C^{4}+63 A^{3} B \,C^{5}+195 A^{3} C^{6}+54 A^{2} B^{7}+354 A^{2} B^{6} C -588 A^{2} B^{5} C^{2}+3423 A^{2} B^{4} C^{3}
$$
$$
+3423 A^{2} B^{3} C^{4}-588 A^{2} B^{2} C^{5}+354 A^{2} B \,C^{6}+54 A^{2} C^{7}+6 A \,B^{8}+108 A \,B^{7} C +354 A \,B^{6} C^{2}+63 A \,B^{5} C^{3}
$$
$$
-378 A \,B^{4} C^{4}+63 A \,B^{3} C^{5}+54 B^{2} C^{7}+6 B \,C^{8}+354 A \,B^{2} C^{6}+108 A B \,C^{7}+6 A \,C^{8}+6 B^{8} C 
$$
$$
+54 B^{7} C^{2}+195 B^{6} C^{3}+357 B^{5} C^{4}+357 B^{4} C^{5}+195 B^{3} C^{6}-63 A^{6} B -63 A^{6} C -441 A^{5} B^{2}-882 A^{5} B C 
$$
$$
-441 A^{5} C^{2}-1008 A^{4} B^{3}-819 A^{4} B^{2} C -819 A^{4} B \,C^{2}-1008 A^{4} C^{3}-1008 A^{3} B^{4}
$$
$$
-3150 A^{3} B^{2} C^{2}-1008 A^{3} C^{4}-441 A^{2} B^{5}-819 A^{2} B^{4} C -3150 A^{2} B^{3} C^{2}-3150 A^{2} B^{2} C^{3}-819 A^{2} B \,C^{4}
$$
$$
-441 A^{2} C^{5}-63 A \,B^{6}-882 A \,B^{5} C -819 A \,B^{4} C^{2}-819 A \,B^{2} C^{4}-882 A B \,C^{5}-63 A \,C^{6}-63 B^{6} C -441 B^{5} C^{2}
$$
$$
-1008 B^{4} C^{3}-1008 B^{3} C^{4}-441 B^{2} C^{5}-63 B \,C^{6}+161 A^{4} B +161 A^{4} C +483 A^{3} B^{2}+791 A^{3} B C 
$$
$$
+483 A^{3} C^{2}+483 A^{2} B^{3}+1260 A^{2} B^{2} C +1260 A^{2} B \,C^{2}+483 A^{2} C^{3}+161 A \,B^{4}+791 A \,B^{3} C 
$$
$$
+1260 A \,B^{2} C^{2}+791 A B \,C^{3}+161 A \,C^{4}+161 B^{4} C +483 B^{3} C^{2}+483 B^{2} C^{3}
$$
$$
+161 B \,C^{4}+16 A^{3}+48 A^{2} B +48 A^{2} C +48 A \,B^{2}+96 A B C +48 A \,C^{2}+16 B^{3}+48 B^{2} C +48 B \,C^{2}+16 C^{3}) \quad .
$$

To see the first $20$ moments (about the mean), enjoy the following output file:

{\tt https://sites.math.rutgers.edu/\~{}zeilberg/tokhniot/oGR3moms4a.txt} \quad .

\vfill\eject

{\bf How we used Symbol-Crunching to Derive these Nice Formulas for the Moments?}

{\bf The Two-Player Case}

Let's first consider the classical case. Let $D(A,B)$ be the random variable

``{\it duration until one of the players goes broke if the initial capitals were $A$ and $B$}'',

and let  $F(A,B)(y)$ be its {\bf probability generating function}.

Obviously, $F(A,B)(y)$ satisfies the following equation:
$$
F(A,B)(y) \, = \, \frac{y}{2} ( F(A+1,B-1)(y)+ F(A-1,B+1)(y)) \quad,
\eqno(1)
$$
subject to the {\it boundary conditions}
$$
F(A,0)(y)=1 \quad, \quad F(0,B)(y)=1 \quad .
$$
Let's expand $F(A,B)(y)$ as a Taylor series around $y=1$,
$$
F(A,B)(y)=\sum_{i=0}^{\infty} f_i(A,B) (y-1)^i \quad.
\eqno(2)
$$
$f_i(A,B)$ is equal to the  {\bf binomial moment}
$$
E \left [ {{D(A,B)} \choose {i}} \right ]  \quad.
$$
Once you have the binomial moments, you can easily compute the {\bf moments}
$$
E[D(A,B)^i]= \sum_{j=0}^{i} j!\,S(i,j)  \,  E \left [ {{D(A,B)} \choose {j}} \right ]  \quad,
$$
where $S(i,j)$ are the Stirling numbers of the second kind.

Once you have the moments, including the most important one, $\mu=E[D(A,B)]$, you can get the {\it moments about the mean} by
$$
E[(D(A,B)-\mu)^i]= \sum_{j=0}^{i} (-\mu)^{i-j}\, {{i} \choose {j}}  E[D(A,B)]^j  \quad.
$$

Plugging-in $(2)$ into $(1)$ and writing, for convenience $y-1=z$, we have
$$
\sum_{i=0}^{\infty} f_i(A,B) z^i 
\, = \, \frac{1+z}{2} \left ( 
\sum_{i=0}^{\infty} f_i(A+1,B-1) z^i 
+
\sum_{i=0}^{\infty} f_i(A-1,B+1) z^i \right ) \quad .
\eqno(3)
$$
Comparing the coefficient of $z^i$, we have the equation
$$
f_i(A,B)= \frac{1}{2} (f_i( A+1,B-1)+f_i(A-1,B+1)) + \frac{1}{2} (f_{i-1} (A+1,B-1)+f_{i-1}( A-1,B+1)) \quad,
$$
implying that
$$
f_i(A,B)- \frac{1}{2} (f_i(A+1,B-1)+f_i(A-1,B+1)) = \frac{1}{2} (f_{i-1}(A+1,B-1)+f_{i-1} (A-1,B+1)) \quad.
$$
Assuming that we already know $f_{i-1}(A,B)$ (of course $f_0(A,B)=1$), we postulate an ansatz for  $f_i(A,B)$ as a {\bf polynomial} in $A,B$ of degree $2i$, using
{\it undetermined coefficients}, letting Maple do the tedious high-school algebra, getting a system of linear equations satisfied by these,
so far undetermined, coefficients. Then we ask Maple to kindly {\tt solve} that system that it generated for us.
To our delight, there was always a unique solution. The uniqueness (in general) is easy by a ``maximum principle".

{\bf The Three-Player (fair) Gambler's Ruin Until One Player Goes Broke}

Everything is analogous. Letting $f_i(A,B,C)$  be the binomial moment of the random variable `duration (until one of the players goes broke)', with initial capitals $A,B,C$, we have
$$
f_i(A,B,C)- 
\frac{1}{6} (
f_i(A+1,B-1,C)+f_i(A-1,B+1,C)
$$
$$
+f_i(A+1,B,C-1)+f_i(A-1,B,C+1) 
$$
$$
+f_i(A,B+1,C-1)+f_i(A,B-1,C+1) )
$$
$$
= \frac{1}{6} 
(f_{i-1}(A+1,B-1,C)+f_{i-1}(A-1,B+1,C)
$$
$$
+f_{i-1}(A+1,B,C-1)+f_{i-1}(A-1,B,C+1) 
$$
$$
+f_{i-1}(A,B+1,C-1)+f_{i-1}(A,B-1,C+1) ) \quad .
$$

Inspired by Engel's formula for the first moment, we took the {\it ansatz} to be a polynomial in $A,B,C$ of degree $2i+1$ divided by $A+B+C$.
Since the denominator is invariant, it is retained, and we can look at the numerator and do the linear algebra with {\it undetermined coefficients},
like we did for the two-player case.

{\bf Duration for a $k$-Player (fair) Gambler's Ruin until Only One Player survives}

The Expectation is easy, as already pointed out by Sheldon Ross [R] (of the popular probability textbooks fame).
If $k$ players have initial capitals $A_1,  \dots, A_k$, the expected duration until $k-1$ of them drop out, is
$$
\sum_{1 \leq i < j \leq k} A_i\, A_j \quad.
$$
Alas, the variance, and higher moments, do not seem to be nice, let alone polynomials (or even rational functions).

{\bf Probability of Winning}

For the classical, two-player, gambler's ruin, famously (and trivially) the probability of winning the whole pot for the player whose
initial capital was $A$ dollars (and the other had $B$ dollars) is $\frac{A}{A+B}$. 
For the {\bf three-player} case, there does not seem to be a nice formula, but 
the fascinating article [DiE] has some very good approximations and asymptotics.
In particular, they proved that with initial  capitals $1,1,N$, the probability is $O(\frac{1}{N^3})$.

Recall that the scaled moments (about the mean) are $m_i/m_2^{i/2}$. 

{\bf The limit of the scaled moments with equal capitals for two  players}

If both players start out with the same capital $a$, then, of course
$$
m_1=E[D]=a^2 \quad, \quad \sigma^2=m_2=\frac{2}{3}a^4 + O(a^3) \quad,
$$
and the limits of the scaled moments (about the mean), from the third to the $20^{th}$, as $a$ goes to infinity,
expressed in decimals, are as follows.
$$
1.9595917942265424786, 8.8285714285714285714, 42.737763893131259770, 255.69610389610389610, 
$$
$$
1775.6935988640577894, 14103.678121878121878, 126010.34642879135344, 1250957.5292367075649, 
$$
$$
13660640.862506747726, 162737708.10778675241, 2100232711.7489454497, 29189813843.019299016, 
$$
$$
434668787296.09441079, 6904215035796.8421121, 116519631313642.16450, 
$$
$$
2082128256599445.4176, 39273259761289810.763, 779763346464591374.25 \quad .
$$
In particular, the {\it limiting skewness} (scaled third-moment about the mean)  is  $1.9595917942265424786$, 
and the {\it limiting kurtosis} (scaled fourth-moment about the mean) is $8.8285714285714285714$, as already computed in [Z].

For the exact values, see the output file

{\tt https://sites.math.rutgers.edu/\~{}zeilberg/tokhniot/oGR3moms5.txt} \quad .

{\bf The Limit of the Scaled moments with equal capitals for three  players}

If all three players start out with the same capital $a$, then, of course
$$
m_1=E[D]=a^2 \quad, \quad \sigma^2=m_2=\frac{1}{2}a^4 + O(a^3) \quad,
$$
and the limits of the scaled moments (about the mean), from the third to the $15^{th}$, as $a$ goes to infinity, in decimals are
$$
1.8384776310850235634, 8.2714285714285714286, 38.719147032686238015, 226.13701298701298701, 
$$
$$
1529.4150671558283928, 11835.953251748251748, 103028.09771971592326, 
$$
$$
996496.58001633815101, 10601990.334510850051, 123051719.84505166456, 
$$
$$
1547212853.2328220329, 20950631863.915761840, 303953682939.86889196 \quad .
$$

In particular, the {\it limiting skewness} is  $1.8384776310850235634$, a little less than the two-player case,
and the {\it limiting kurtosis} is $ 8.2714285714285714286$, also a little less than its counterpart.
For the exact values, see the output file

{\tt https://sites.math.rutgers.edu/\~{}zeilberg/tokhniot/oGR3moms6.txt} \quad .

{\bf More than Three-Players?}

With the exception of Ross's formula (mentioned above), that only handles the expectation, for the $k$-player game until only one player survives,
our approach does not seem to work with more than three players.
Of course, perhaps there is yet-another ansatz, and perhaps one can get, by {\it symbol-crunching},
at least asymptotic results, but this remains to be seen. Again we strongly recommend the very interesting Diaconis-Ethier article [DiE], of whose
existence we learned thanks to Jean-Paul Delahaye's engaging column.

{\bf References}

[A] Jiri Andrei, {\it ``Mathematics and Chance''}, Wiley, 2001.

[BLT] F.T. Bruss, G. Louchard, and J.W. Turner, On the $N$-tower problem and related problems. Adv. Appl. Probab. {\bf 35} (2003), 278-294.

[De] Jean-Paul Delahaye, {\it Ruiner le casino ou se ruiner}, Pour la Science {\bf 550} (August 2023), 80-85. \hfill \break 
{\tt  https://medias.pourlascience.fr/api/v1/files/64b5166b903f8e26c07a14f5?alt=file}

[DiE] Persi Diaconis and  Stewart N. Ethier, {\it Gambler’s Ruin and the ICM}, Statist. Sci. {\bf 3}, 289-305, August 2022. 
{\tt https://arxiv.org/abs/2011.07610} \quad .

[E] Arthur Engel, {\it The computer solves the three tower problem}, Amer. Math. Monthly {\bf 100} (1993), 62-64.

[F] William Feller, ``{\it An Introduction to Probability
Theory and Its Application}'', volume 1, three editions. John Wiley and sons. First edition: 1950. Second edition:1957. Third edition: 1968.
Volume 2: two editions.

[R] Sheldon M. Ross, {\it A simple solution to a multiple player gambler's ruin problem}, Amer. Math. Monthly {\bf 116} (2009), 77-81.

[S] David Stirzaker, {\it Three-handed gambler's ruin}, Adv. Appl. Probab. {\bf 38} (2006), 284-286.

[Z] Doron Zeilberger, {\it Symbol-crunching with the gambler’s ruin problem},  in: ``Tapas in Experimental Mathematics", Contemporary Mathematics {\bf 457} (2008), 285-292, (Tewodros Amdeberhan and Victor Moll, eds.) \hfill\break
{\tt https://sites.math.rutgers.edu/\~{}zeilberg/mamarim/mamarimhtml/ruin.html} \quad .

\bigskip
\hrule
\bigskip

Shalosh B. Ekhad, c/o D. Zeilberger, Department of Mathematics, Rutgers University (New Brunswick), Hill Center-Busch Campus, 110 Frelinghuysen
Rd., Piscataway, NJ 08854-8019, USA. \hfill\break
Email: {\tt ShaloshBEkhad at gmail  dot com}   \quad .
\smallskip

Doron Zeilberger, Department of Mathematics, Rutgers University (New Brunswick), Hill Center-Busch Campus, 110 Frelinghuysen
Rd., Piscataway, NJ 08854-8019, USA. \hfill\break
Email: {\tt DoronZeil at gmail  dot com}   \quad .
\bigskip
Sept. 15, 2023.

\end